\documentclass[showpacs,showkeys,preprintnumbers,amsmath,amssymb]{revtex4}

\usepackage{graphicx}
\usepackage{dcolumn}
\usepackage{bm}

\begin{document}


\title{Fibonacci-like sequences and shift spaces in symbolic dynamics}

\author{R. Tonelli}
 \altaffiliation[Also at ]{SLACS-CNR-INFM and INFN at Physics Department, Cagliari University.}
 \email{roberto.tonelli@dsf.unica.it, tonelli@eecs.berkeley.edu}
\affiliation{%
CNR-SLACS-INFM and INFN, Physics Dept., University of Cagliari - Italy
}%


\date{\today}

\begin{abstract}
We afford the problem of counting the blocks of a given length 
made with symbols drawn from an alphabet and relate this number 
to Fibonacci-like recurrent relations. 
The recurrence polynomia allows to calculate the limit ratio of two adjacent terms 
and this information is used to determine the topological entropy of a discrete numbers of associate shift spaces.
We then describe a scheme to build a shift space with a pre-selected entropy. 
 
%
%
%
\end{abstract}

\maketitle

\section{Introduction}
Modern technology and, more in general, modern society, relay on computer systems 
in which information is processed and represented in form of binary strings.
The problem of transmission and storage of huge amount of data in the form 
of binary series has become one of the main 
interest of information science and modern mathematics.
Usually these binary series are subjected to some logical constraint
so that not all the binary sequences may be allowed and 
there may be sets of forbidden 'words'.
Practical examples of this situations are the memorization of information 
onto magnetic supports such hard disks, where these sequences are stored 
in a series of tracks which must have a 'preamble' or a 'end' identified 
by a well defined sequence of bits; sometimes it is convenient 
to fill the space among two strings with a succession of zeroes; 
error detection and correction codes, both in the problem of storage and 
transmission of information, require some sequences of bits to be forbidden. 
Another main field of research in which bits sequences require 
to be subjected to some constraint are the problems connected 
to the security of data stored and transmitted. When a computer system 
is protected through a password or when the data transmission is performed 
through cryptography it may be of relevance to know 
how many different sequences, among the bit sequences of a given length and subjected 
to constraints, may be generated. 
Another interesting discipline in which this kind of problems 
are of relevance is the field of dynamical systems. 
Symbolic dynamics allows the representation of the time evolution of 
a dynamical system through a discrete succession of two or more symbols 
in which all the essential information about the dynamical 
system, and in general its physical properties, are preserved. 
One example is the concept of entropy that, for a given symbolic 
representation of a dynamical system may be directly related 
to the rate of growth of the number of different sequences with the length 
of the sequence.

In this paper we discuss some of these aspects using the language of {\em shift spaces}.
In particular we show 
how many of these problems may be expressed and solved in terms of 
Fibonacci-like recurrent relations in a simple, intuitive and amenable way.
Finally we relate the limit ratios of such sequences to the topological entropy 
of the corresponding shift space and show how can be created a shift space with a predetermined entropy.

%
%
%
%

%

\section {Shift Spaces and Forbidden Blocks}
All of the above mentioned problems may be better expressed in a rigorous 
mathematical language through the concept of {\em shift space} \cite{LM}. 
Even if the main interest of this paper regards binary sequences, some results 
are immediately extendable to sequences of more than two symbols, thus we start 
our discussion recalling the concept of shift space over a finite alphabet of $k$ symbols.

Given a finite alphabet ${\cal A}$ = $\{a_{0}, a_{1}, \dots, a_{k-1}\}$ of $k$ symbols 
and the set of integers $\mathbb{Z}$ the {\em full shift} is defined as the set of bi-infinite sequences 
of ${\cal A}^{\mathbb{Z}}$ defined by:
\begin{eqnarray}
{\cal A}^{\mathbb{Z}} = \{ x = ( x_{i} )_{i \in \mathbb{Z}} : x_{i} \in A ~for ~all ~i \in \mathbb{Z} \}
\end{eqnarray}
namely {\it the set of all functions from $\mathbb{Z}$ to A}.
On this set a {\em left shift} is naturally defined as the function $\sigma$ defined by 
\begin{eqnarray}
\sigma : {\cal A}^{\mathbb{Z}} \rightarrow {\cal A}^{\mathbb{Z}}  ~:~ 
\sigma (x)_{i} = x_{i+1}~~  ~for ~all ~i~ \in \mathbb{Z}
\end{eqnarray}
A finite sequence of (length-$m$) symbols drawn from ${\cal A}$
is called a {\em block over} ${\cal A}$:
\begin{equation}
m-block = x_i \dots x_{i+m},  
\end{equation}
and its length is indicated by the norm $\vert x_i \dots x_{i+m} \vert = m$.
An m-block made by repetition of a single symbol is indicated $a^m$ where $a \in {\cal A}$ is repeated $m$-times. 
A {\em shift space {\cal X}} is a subset of the full shift in which not all 
the sequences of $k$ symbols are allowed.
A {\em shift space} may be identified or defined in two main and equivalent ways.
One can give the set ${\cal F}_X$ of all forbidden blocks 
in the sequences of symbols, or one can give what is called the {\em language}
of the shift space, namely the collection $B_{n}(X)$ of all the allowed words of any length $n$.
Here are some examples.

X is the set of all binary bi-infinite sequences where the block (11) is not allowed. 
This is called the {\em golden mean shift}. 

X is the set of all binary sequences such that each pair of 'ones' is separated 
by an even number of 'zeroes', thus the collection of forbidden blocks takes the form:
${\cal F} = \{ 10^{2n+1}1 : n \geq 0 \}$.

Thus in any shift space that is not the full shift not all the sequences
of symbols will appear. 
Consider as the simplest case an alphabet of two symbols that we identify 
through the two bits 'one' and 'zero'.
If we refer to the arguments introduced above we see how any of the 
problems discussed may be related to the identification of the 
appropriate shift space.
Furthermore it is easy to see that in the case of dynamical 
systems forbidden and allowed blocks correspond to forbidden 
or allowed transitions among different physical states and to 
forbidden and allowed physical behaviors.
In particular one can calculate the {\em topological entropy} \cite{topological} of a 
shift space associated to a dynamical system by counting 
the different allowed words of a given length and examining 
how this number grows when the length tends to infinity \cite{LM}:
\begin{equation}
h = \lim_{n \rightarrow \infty} \frac{log(B_{n}(X))}{n}
\label{eq:entropy}
\end{equation}
In the next section we relate these numbers to different series 
of Fibonacci-like sequences. 

\section{The number of allowed words of length $n$ for different shift spaces}
\subsection{Golden mean shift and Fibonacci sequence}
The {\em Golden mean shift} is a shift space defined by only one 
forbidden block: 
\begin{equation}
{\cal F} = \{11\}
\end{equation}
Examining all the sequences allowed one may check that they 
consist only of isolated ones separated apart by at least one zero. 
One may think of this shift space as of all the sequences 
of bits stored onto an hard disk in which zeroes are need for timing purposes, or 
of all the transmitted sequences of bits with Humming distance 1, or 
of one dynamical system in which the trajectory must escape 
immediately form one of the two possible states in which it may exist and so on.

Thus it is interesting, both for practical purposes and from 
a theoretical point of view, to calculate the number \#$(n)$ of different strings 
of length $n$ of such form.
Here we show first by induction and next by a constructive scheme 
that this number is exactly the $n^{th}$ Fibonacci number: \#$(n)$ = $F_n$ 
where the relation $F_n = F_{n-1} + F_{n-2}$ holds.

{\bf Lemma 1}
Let $\#(n)$ be the number of different blocks of length $n$ for the golden mean shift. 
Then \#$(n)$ = $F_n$ where $F_n$ is the $n^{th}$ Fibonacci number.

{\em Proof:}
For $n = 1$ is \#$(1)$ = 2 and for $n = 2$ is \#$(2)$ = 3, thus we 
deal with Fibonacci numbers with initial conditions $F_0 = 1$ and 
$F_1 = 2$ (alternatively, we may choose to set \#$(n)$ = $F_{n+1}$ with 
initial conditions $F_0 = 1$ and 
$F_1 = 1$). 

The first step is to check that \#$(3)$ = \#$(2)$ + \#$(1)$ = 5, and this is 
done immediately.
Next, supposing that \#$(n)$ = \#$(n-1)$ + \#$(n-2)$ is valid, we show that 
\#$(n+1)$ = \#$(n)$ + \#$(n-1)$ must be valid too.
In fact, given a string of $n$ bits 
we may construct all the allowed strings of length $n+1$ 
of the same type by adding the bit 'zero' 
or the bit 'one' at the end of each length-$n$ allowed string. 
Adding a 'zero' presents no problems. It may be added to any of the 
\#$(n)$ strings obtaining \#$(n)$ strings of length $n+1$ terminating by 'zero'.
The bit 'one' may be added only to strings of length $n$ terminating by 'zero', since 
the couple (the block) $(11)$ is not allowed.
Thus adding 'one' at the end we obtain a number of strings 
of length $n+1$ terminating by the two block word $(01)$.

Now we may think of these strings as compound strings of initial $n-1$ bits and of two 
final bits equal to $(01)$. The $n-1$ bits part, since it is followed by 'zero', 
has no other constraint than the forbidden block $(11)$ and thus 
there are \#$(n-1)$ of such compound strings. Thus the "adding one" process 
provides exactly \#$(n-1)$ strings of length $n+1$ ending by the word $(01)$.
The total number of allowed strings of length $n+1$ is given by those obtained adding 
'zero' plus those obtained adding 'one', namely we get:
\begin{equation}
\#(n+1) = \#(n) + \#(n-1)
\end{equation}
after the hypothesis \#$(n)$ = \#$(n-1)$ + \#$(n-2)$. 

\begin{table}
\caption{}
\begin{ruledtabular}
\begin{tabular}{cccc}
\#(1) = 2 & \#(2) = 3 & \#(3) = 5 & \#(4) = 8\\
\hline
\\
0 & 00 & 00~~0 & 000~~0\\
1 & 01 & 01~~0 & 010~~0\\
  & 10 & 10~~0 & 100~~0\\
  &    & 0~~01 & 001~~0\\
  &    & 1~~01 & 101~~0\\
  &    &      & 00~~01\\
  &    &      & 01~~01\\
  &    &      & 10~~01\\

\end{tabular}
\end{ruledtabular}
\label{tab-fibonacci}
\end{table}

All the reasoning will be even clearer when constructing from scratch 
the recursive relation.
A look at table \ref{tab-fibonacci} shows a constructive demonstration.
In the columns from one to four are illustrated all the different configurations of an $n$-block, $n = 1, 2, 3, 4$
for the {\em golden mean shift}. The construction process starts with columns one and two, where 
the cases with $n = 1,2$ are reported. 
Next we build all the $3$-blocks starting from all the allowed blocks 
with $n = 2$. 
We may obtain such blocks adding the bit 'one' or the bit 'zero' to all the \#$(2) = 3$ strings with $n = 2$.
The bit 'zero' may be added at the end of each string, giving the first \#$(2) = 3$ blocks with $n = 3$ ending by 'zero'.
The bit one, due to the exclusion of the block $(11)$, may be added only to the $2$-blocks ending 
by 'zero'. Thus any string obtained in this last way will have length $n = 3$ and will end with the block $(01)$. 
Then we may think of all such strings as made by a first part, of length $3 - \vert (01)\vert = 3 -2 = 1$, 
where all the possible configurations allowed for the {\em golden mean shift} are present, and the last part 
made by the block $(01)$. The number obtained in this way is of course \#$(3 - \vert (01)\vert)$ = 2. 
The total number of newly generated strings will thus be \#$(2)$ + \#$(1)$. 
The same reasoning is applied when we want to build all the allowed $4$-blocks starting 
from all the allowed $3$-blocks. In table \ref{tab-fibonacci} the strings are 
represented leaving a space in the sequences in order to easily identify 
the strings obtained adding 'zero' and those obtained adding 'one'.

\subsection{$T_2$-Shift and Lam\'e numbers.}
In the {\em golden mean} shift each 'one' is isolated and separated by at leas one 'zero' 
from others. Here we examine a shift space in which each 'one' is again isolated 
but separated by at least two 'zeroes' from others. We call it $T_2$-shift \cite{nota} and it is 
described through the following forbidden blocks:
\begin{equation}
{\cal F}_{T_2} = \{011, 101, 110, 111\}  \\
\end{equation}

Again, these numbers may represent a scheme for storing bit sequences 
where at least two 'zeroes' must separate each single 'one', like  
all the transmitted sequences of bits with Humming distance 2, 
or one dynamical system in which the trajectory must escape 
immediately form one of the two possible states in which it may exists 
and presents a latent time of at least two before returning into the same state, like 
a quantum system with a discrete set of states each one with its own decay time and so on.

We calculate the number \#$(n)$ of different strings 
of length $n$ of such form. 
This number is given by a recurrence series similar to the Fibonacci numbers, where 
the relationship is: \#$(n)$ = $F_n$
with $F_n = F_{n-1} + F_{n-3}$. These numbers constitute a series known as Lam\'e series \cite{Lame}
that is: $1, 2, 3, 4, 6, 9, 13, 19, 28. \dots $.

\begin{table}

\caption{}
\begin{ruledtabular}
\begin{tabular}{ccccc}
\#(1) = 2 & \#(2) = 3 & \#(3) = 4 & \#(4) = 6 & \#(5) = 9\\
\hline
\\
0 & 00 & 000 & 000~~0 & 0000~~0\\
1 & 01 & 010 & 010~~0 & 0100~~0\\
  & 10 & 100 & 100~~0 & 1000~~0\\
  &    & 001 & 001~~0 & 0010~~0\\
  &    &     & 0~~001 & 0001~~0\\
  &    &     & 1~~001 & 1001~~0\\
  &    &     &        & 00~~001\\
  &    &     &        & 01~~001\\
  &    &     &        & 10~~001\\

\end{tabular}
\end{ruledtabular}
\label{tab-lame}
\end{table}
We first construct a table \ref{tab-lame} for this case and then provide a demonstration by induction. 
We start writing the first three columns in which the allowed strings with $n$ ranging 
from one to three are reported. Next we want to obtain all the allowed strings with 
$n = 4$ and the constrain that each appearing 'one' must be separated from others 
by at least two 'zeroes'. 
Starting from all the $3$-blocks we may add the bits 1 or 0 at the end. 
The adding 'zero' step is allowed for any $3$-block and gives \#$(3) = 4$ blocks of 
length $n = 4$ satisfying the constraints. 

The digit one may be added only to the $3$-blocks
ending with at least two 'zeroes', namely with the block $(00)$, thus providing 
$4$-blocks ending by the block $(001)$ where the preceding part is 
of length $(4 - \vert (001)\vert) = 4 -3 = 1$, in which each allowed configuration must be present.
 
The reasoning is perfectly similar to the previous section. 
Thus, when adding the digit 'one', we obtain a number \#($4 - \vert (001)\vert) = \#(4 - 3) = \#(1)$ of different 
allowed strings. The total is \#$(4)$ = \#$(3)$ + \#$(1)$. 
The same happens when we start from $4$-blocks to build allowed $5$-blocks, 
adding the digits 'zero' and 'one'. Adding zero gives the new \#$(4)$ strings ending by 'zero' 
and adding 'one' corresponds to add the block $(001)$ to all the strings of length 
$5 - 3 = 2$. The new number is \#$(5)$ = \#$(4)$ + \#$(2)$ and the recurrence relationship
is \#$(n)$ = \#$(n-1)$ + \#$(n-3)$. Thus we have the following \\
{\bf Lemma 2} Let $\#(n)$ be the number of different blocks of length $n$ for the $T_2$ shift. 
Then \#$(n)$ = $F_n$ where $F_n$ is the $n^{th}$ Lam\'e number and the recursive relationship 
\#$(n)$ = \#$(n-1)$ + \#$(n-3)$ holds.

{\em Proof:}
It is \#$(1) = 2$, \#$(2) = 3$, \#$(3) = 4$, thus we are using Lam\'e numbers with these initial conditions.
We already checked that \#$(5)$ = \#$(4)$ + \#$(2)$. 
Next, given \#$(n)$ = \#$(n-1)$ + \#$(n-3)$, to obtain the number of allowed strings 
of length $n + 1$ we have to add a digit to the \#$(n)$ allowed strings with $n$ digits. 
If the digit is $0$ there are no problems and we obtain new \#$(n)$ strings of length $n + 1$. 
Adding the digit $1$ we must eliminate all the strings ending not by $(00)$. 
Then all these blocks must terminate by the three digits block $(001)$ and may be thought 
as made of this last three digits and by a first $n +1 - 3$ digits part in which 
all the allowed \#$(n-2)$ configurations which satisfy the given constraints are allowed. 
Thus, under the hypothesis made, we obtain \#$(n+1)$ = \#$(n)$ + \#$(n-2)$ and 
the recurrence relationship holds for any $n$.

\subsection{$T_m$-Shifts with two symbols.}
Let denote with $T_m$ the shift spaces where ${\cal A}$ = \{0,1\} with the constrain 
that the allowed sequences are made of isolated 'ones' separated by at least $m$ 
'zeroes'.
It is straightforward to generalize to such shift spaces 
the recurrence relationship as:
\begin{equation}
\#(n) = \#(n-1) + \#(n-m)
\end{equation}

{\bf Lemma 3} Let $\#(n)$ be the number of different blocks of length $n$ for the $T_m$ shifts. 
Then \#$(n)$ = $F_n$ where $F_n$ is the $n^{th}$ number of the recursive relationship 
$F_n = F_{n-1} + F_{n-m}$ with opportune initial conditions.

{\em Proof:}
To obtain all the sequences of length $n + 1$ from the allowed $n$-blocks we may add 'zero', getting 
new \#$(n)$ allowed strings, or we may add 'one'.
To add 'one' we must take into account that it must be added 
only to those strings ending by $m$ 'zeroes', obtaining $(n + 1)$-blocks 
ending with the block $0^m1$ that may be added to 
any one of the \#$(n-m)$ allowed blocks of length $n - m$, and this gives new \#$(n-m)$ strings.

It is worth noting that in the limit $n \rightarrow \infty$ all of these 
recurrence series may be solved using the correspondent recurrence equation \cite{Bat}
to obtain the limit ratio $F_n / F_{n-1} = K$ so that the $n^{th}$ terms 
may be expressed as a limit $F_n \simeq F_0 \times K^n$. 

For the Fibonacci series this limit is the {\em golden section} $\phi = \frac{1 + \sqrt{5}}{2}$ 
that may be obtained from the corresponding polynomial recurrence equation:
\begin{eqnarray}
\label{root16}
\lambda^{n} = \lambda^{n-1} + \lambda^{n-2} \Rightarrow \\ \nonumber
\Rightarrow  \lambda^{2} = \lambda + 1
\end{eqnarray}
as the greatest root. 

For the other recurrence series it holds a similar rule so that the limit may be 
recovered as the greatest root of the corresponding polynomial recurrence equation:
\begin{equation}
\lambda^{m+1} = \lambda^{m} + 1
\end{equation}

\section{More than two symbols and Different constraints.}
Next we briefly examine some cases in which the alphabet contains more than 
two symbols or has different constraints. 
\subsection{k-symbols, same kind of constraints}


Let ${\cal A}$ contain $k$ symbols, that we assume to be $\{0, 1, \dots, k-1\}$. 
We set as constraint the condition that each symbol, except 'zero', must appear isolated 
and separated by at least $m$ 'zeroes' from symbols other than 'zero'. 
We indicate these shift spaces $T_{m,k}$.

{\bf Lemma 4}
To these shift spaces are associated recurrence series (with opportune initial conditions) which give 
the number of different allowed configurations (~$\#(n)$~) for a string of length $n$ according to 
the relationship:
\begin{equation}
\#(n) = \#(n-1) + (k-1)*\#(n-m-1)
\label{eq:n}
\end{equation}
{\em Proof:} 
To create the ($n+1$)-blocks we add 
one of the $k$ symbols to the blocks of length $n$. 
'Zero' may be added to any block giving $\#(n)$ new ($n+1$)-blocks. 
For each of the $k-1$ symbols different from 'zero' we obtain strings of length $n+1$ made of two parts:
the first part has length $n-m$ and is one of the $\#(n-m)$ allowed strings of that length; the last part 
contains one of the $k-1$ blocks of the form $0^{m}a$, where $a$ is one of the $k -1$ symbols different from 'zero'. 
Thus there are a total of $k-1$ possibilities of adding these final blocks in queue of
any allowed $(n-m)$-block.
The total provides the number of length-$(n+1)$ allowed blocks as $\#(n+1) = \#(n) + (k-1)*\#(n-m)$ 
that is the same as eq.~(\ref{eq:n}).
For $k = 2$ and $m = 1$ one recovers the Fibonacci series and for $k = 2$ and $m = 2$
the Lam\'e series.

\begin{table}
\caption{}
\begin{ruledtabular}
\begin{tabular}{ccccc}
m,k = 1,1 & 1,2 & 1,3 	& 1,21		& 2,5	\\
\hline
\\
2 	& 2 	& 3 	& 21 		& 5	\\
3 	& 3 	& 5 	& 41 		& 9	\\
5 	& 4 	& 11 	& 461 		& 13	\\
8 	& 6  	& 21 	& 1281 		& 33	\\
13 	& 9   	& 43  	& 10501 	& 69	\\
21 	& 13  	& 85  	& 36121		& 121	\\
34 	& 19   	& 171  	& 246141	& 250	\\
55 	& 28   	& 341 	& 968561	& 526	\\
\vdots  & \vdots &\vdots &\vdots &\vdots\\
\hline
$F_n/F_{n-1} \rightarrow \phi$ & 2 & 2 & 5 & 2 \\ 

\end{tabular}
\end{ruledtabular}
\label{tab-series}
\end{table}

All the recurrence series of such form have a limit ratio among two adjacent terms 
when $n \rightarrow \infty$ and this limit is the largest root (indicated by $\lambda_0$ from now on) 
of the polynomia associated to the recurrence series, as can be obtained with the same technique used for the 
Fibonacci series:
\begin{eqnarray}
\lambda^n = \lambda^{n-1} + (k-1)\lambda^{n-m-1}   \Rightarrow \\ \nonumber
\Rightarrow  \lambda^{m+1} = \lambda^m + (k - 1)
\label{eq:final}
\end{eqnarray}
This allows to select the values for the couple $m, k$ such that 
the largest root has a prefixed value. 
Consider for example the cases 
(a) $m = 1, k = 3$, (b) $m = 1, k = 21$, (c) $m = 2, k = 5$.
In table \ref{tab-series} are listed these and some other series for various initial conditions.

Let consider (a). The limit ratio among two adjacent terms tends to the value 2. 
In fact the equation~(\ref{eq:final}) becomes $\lambda^{2} = \lambda^1 + (3 - 1)$
with largest root $\lambda_0 = 2$. 

For example one can use $m = 2$ and set $\lambda^{2} = \lambda^1 + (k - 1)$. Then 
choosing a given value for $\lambda_0$ one can solve for $k$ such that the limit ratio is that 
value of $\lambda_0$ (with $k$ natural number). 
The other examples are (b), such that the limit ratio be 5, and (c), such that 
the limit ratio be again 2, but with different $m$. The formula (\ref{eq:final}) leaves 
ample freedom for fixing the desired $m$ and $\lambda$ and solve for the correspondent $k$. 

\subsection{Three symbols, different kind of constraints}
Next we fix ${\cal A}$ = \{0, 1, 2\} and we allow for the symbols 'one' and 
'two' to be next to each other, such that ${\cal F}$ = \{11, 22\}. 
In this case, given all the strings of length $n$, we may obtain the strings of 
length $n+1$ adding 'zero', 'one' or 'two' in queue. 
Adding 'zero' provides new $\#(n)$ strings. 
The situation for the digits 'one' and 'two' is symmetric, thus we need only to consider 
one of the two. 
Let consider the addition of the digit 'one'. It may be added to all the strings ending by 'zero' 
or ending by 'two', to get the final blocks (01) or (21). 
The final block (01) constitutes the last part of a length-$(n+1)$ string 
in which the first part, of $n-1$ digits, has no other constraints than those given by ${\cal F}$. 
There are $\#(n-1)$ allowed such configurations. The symmetric holds when considering 'two' as last digit, thus there are another 
$\#(n-1)$ allowed such configurations.
If we consider the block (21) it may be the final part only of those ($n+1$)-blocks in which 
the first $n-1$ block does not end by 'two'. If it ends by 'zero', we get a final block of three digits 
(021) and a first part of length $n-2$ where there are $\#(n-2)$ different allowed configurations. 
The symmetric holds again for the digit 'two', giving a three block (012) which provides itself 
other $\#(n-2)$ allowed configurations. 
If it ends by 'one', we get a final block of three digits (121) that may be added in queue only 
to blocks ending not by 'one'. The same symmetrically for the three digits (212). 
We obtain two branches of growing final block with alternating sequence ...121212... and vice-versa, 
from which at each back-step depart $2*\#(n-i)$ allowed configurations, for the index $i = 1, 2, \dots, n-1$. 
Then we get $\#(n+1) = \#(n) + 2\#(n-1) + 2\#(n-2) + \dots 2\#(1) + 2 + 2$ 
where the last two numbers indicate respectively the strings $0121212....$, $0212121....$ and 
the two complementary strings $121212....$, $212121....$, all blocks of $n+1$ digits.  
Then we have:
\begin{equation}
\#(n) = \#(n-1) + 2\sum_{i = 2}^{n -1} \#(n-i) + 4.
\end{equation}


\section{Topological Entropy for $T_{m,k}$ shift spaces}
From the previous sections and from the equation (\ref{eq:entropy})
it is clear how one can proceed to select a $T_{m,k}$ shift space which possesses 
a desired value for the entropy. 
For the examined cases it is $B_n(X) = \#(n)$. 
To evaluate the limit for $n \rightarrow \infty$ it is convenient 
to express $\#(n)$ through the limit ratio $F_n/F_{n-1}$, so that 
\begin{eqnarray}
\#(n) \simeq \lambda_0^{n-1} \times \#(1) \\
\\ \nonumber
h = \lim_{n \rightarrow \infty} \frac{log(\lambda_0^{n-1} \times \#(1))}{n} = log(\lambda_0)
\end{eqnarray}
Using the previous results, we take for example $m = 1$ and let $k$ be arbitrary. 
The resulting equation for the largest root $\lambda_0$ is:
\begin{eqnarray}
\lambda^2 = \lambda + (k-1) \Rightarrow \\
\Rightarrow
\lambda_0 = \frac{1 + \sqrt{4k -3}}{2}
\end{eqnarray}
and we get a discrete set of numerable topological entropies:
\begin{equation}
h_{T_{2,k}} = log( \frac{1 + \sqrt{4k -3}}{2} )
\end{equation}
The freedom we have in choosing the two integers $k$ and $m$ allows 
for a vast variety of $T_{m,k}$ shift spaces with a pre-selected associated entropy.

%
%
%
%

\section{Conclusions}
We showed how the number of different allowed blocks of fixed length 
for a large class of shift spaces ($T_{m,k}$) may be expressed in terms of 
recurrence relations, obtaining various series describing the growth of this 
number with the length of the blocks. 
This information has been used to calculate different topological entropies
from the limit ratio of two adjacent terms of such series. 
A scheme to build a desired value for the entropy has been also illustrated. 
It must be noted that the recurrence relations we obtained for particular shift spaces, and 
used to pre-select some values for $h_{T_{m,k}}$, may have a larger applicability
to problems arising in different contexts. 
In fact it is easy to check that they may be equivalently used to count 
the number of different one-dimensional paths (if $\vert {\cal A} \vert = 2$) 
where one must step at least $m$ times to the right before taking any step to 
the left for only once. 
Finally we believe that all these results may be applied when devising complex algorithms for information 
storage and encryption, since the entropy is directly related to the complexity 
of the system.











\section*{Acknowledgments}
This work has been partially supported by the R.A.S (Regione Autonoma Sardegna), under the project M\&B-T12

\newpage 
\bibliography{apssamp2math}

\begin{thebibliography}{5}
\bibitem[{not()}]{nota}
\eprint{The notation $T_{m,k}$ comes from the author name, since to our
  knowledge no such general recurrence series and shift spaces have been used.}
\end{thebibliography}

\end{document}